\documentclass[a4paper,10pt]{article}
\usepackage{amsfonts}
\usepackage{amstext}
\usepackage{amsthm}
\usepackage{amsmath}
\usepackage{amssymb}
\usepackage{latexsym}
\usepackage[latin1]{inputenc}
\newcommand{\bl}{\hfill\rule{2mm}{2mm}}
\newcommand{\R}{\Bbb{R}}

\newtheorem{mainthm}{Theorem}

\newcommand{\n}{\noindent}

\begin{document}

\title{Sharp constants in Riemannian $L^p$-Gagliardo-Nirenberg inequalities
 }
\author{\textbf{Jurandir Ceccon \footnote{\textit{E-mail addresses}:
ceccon@ufpr.br (J. Ceccon)}\hspace{0,2cm}; \textbf{Carlos E. Dur\'{a}n \footnote{\textit{E-mail addresses}: cduran@ufpr.br (C. Dur\'{a}n)}}}\\
{\small\it Departamento de
Matem\'{a}tica, Universidade Federal do Paran\'{a},}\\
{\small\it Caixa Postal 019081, 81531-990, Curitiba, PR,
Brazil}}\maketitle

\markboth{abstract}{abstract}
\addcontentsline{toc}{chapter}{abstract}

\hrule \vspace{0,2cm}

\n {\bf Abstract}

Let $(M,g)$ be a smooth compact Riemannian manifold of dimension
$n \geq 2$, $1 < p < n$ and $1 \leq q < r < p^\ast = \frac{np}{n-p}$ be real parameters.
This paper concerns to the validity of the optimal
Gagliardo-Nirenberg inequality

\[
\left( \int_M |u|^r\; dv_g \right)^{\frac{\tau}{r \theta}} \leq
\left( A_{opt} \left(\int_M |\nabla_g u|^p\;
dv_g\right)^{\frac{\tau}{p}} + B_{opt} \left(\int_M |u|^p\;
dv_g\right)^{\frac{\tau}{p}} \right) \left( \int_M |u|^q\; dv_g
\right)^{\frac{\tau(1 - \theta)}{\theta q}} \; .
\]

This kind of inequality is studied in Chen and Sun (Nonlinear Analysis 72 (2010), pp. 3159-3172)
where the authors established its validity when $2 < p < r < p^\ast$
and (implicitly) $\tau = 1$. Here we solve the case $p \geq r$ and introduce
one more parameter $1 \leq \tau \leq \min\{p,2\}$.  Moreover, we prove the existence of extremal
function for the optimal inequality above.

\n
\vspace{0,5cm} \hrule\vspace{0.2cm}

\section{Introduction}

Inequalities of Gagliardo-Nirenberg type (\cite{Ga} and
\cite{Ni}) contain by varying the parameters some classical
inequalities, such as the  Moser \cite{Mo} and Nash \cite{Na} inequalites.
Moreover, by taking limits of the parameters, we get logarithmic \cite{DPDo} and
Sobolev \cite{S} inequalities. 

Optimal inequalities of Gagliardo-Nirenberg type have been extensively studied, both in the
Euclidean and Riemannian contexts; e.g, \cite{Ma}, \cite{B1},
\cite{CaLo}, \cite{CoVi}, \cite{DelDo} for the Euclidean case and
\cite{Au1}, \cite{AuLi}, \cite{Ba}, \cite{Bro}, \cite{CMMZ}, \cite{CMJD}, \cite{CS}, \cite{D2}, \cite{Hu} for
Riemannian manifolds.

The {\em optimal} cases of Euclidean Gagliardo-Nirenberg
inequalities are used, for example, for finding sharp criteria for the
global existence for nonlinear Schr\"odinger equations (see \cite{CG}
or \cite{ShuZhang}), and optimal decay rate of the intermediate
asymptotics of solutions to nonlinear diffusion equations (see
\cite{DelDo}). Recently, the Riemannian Gagliardo-Nirenberg optimal constants
studied in \cite{CMMZ} where applied by \cite{KM} to obtain
global existence theorems for Zakharov system in $\mathbb{T}^2$.
A particularly important family of applications of optimal Gagliardo-Nirenberg inequalities is the transition 
to optimal {\em Entropy inequalities},  in the spirit of 
\cite{ceccon-montenegro-logarithmic, DPDo}.

Denote by $D^{p,q}(\R^n)$ the completion of $C_0^\infty(\R^n)$
under the norm

\[
||u||_{D^{p,q}(\R^n)} = \left( \int_{\R^n} |\nabla u|^p\; dx
\right)^{\frac{1}{p}} + \left( \int_{\R^n} |u|^q\; dx
\right)^{\frac{1}{q}}\, .
\]

\n The Euclidean Gagliardo-Nirenberg inequality states that there exists $A > 0$,
such that for any function $u \in D^{p,q}(\R^n)$,

\begin{gather}\label{dgne}
\left( \int_{\R^n} |u|^r\; dx \right)^{\frac{p}{r \theta}} \leq A
\left( \int_{\R^n} |\nabla u|^p\; dx \right) \left( \int_{\R^n}
|u|^q\; dx \right)^{\frac{p(1 - \theta)}{\theta q}}\, ,
\tag{$GN_E(A)$}
\end{gather}

\n where $1 < p < n$, $1 \leq q < r < p^* = \frac{np}{n - p}$ and
$\theta = \frac{np(r - q)}{r(q(p - n) + np)} \in (0,1)$ is the
interpolation parameter. Define

\[
A(p,q,r,n)^{-1} = \inf_{u \in D^{p,q}(\R^n)} \{ ||\nabla
u||_{L^p(\R^n)}^p ||u||_{L^q(\R^n)}^{\frac{p(1 -
\theta)}{\theta}}; ||u||_{L^r(\R^n)} = 1\} \; .
\]

\n Note that this constant is well defined since the right hand side infimum is a positive number by $GN_E(A)$. Inequality
$GN_E(A(p,q,r,n))$ is called {\em optimal Euclidean
Gagliardo-Nirenberg inequality} and the constant $A(p,q,r,n)$ is the {\em best
constant} in this inequality. A non-zero function realizing equality in $GN_E(A(p,q,r,n))$
is said to be an {\em extremal function}. The existence of such function is
established in this case by using standard classical methods of Calculus of Variations.

We now consider the Riemannian case. Let $(M,g)$  be a smooth, compact Riemannian manifold without boundary of
dimension $n \geq 2$. Using standard arguments, we obtain a Riemannian version of the Euclidean inequality $GN_E(A)$:
there exists positive constants $C,D$ such that for all $u$ in the Sobolev space $H^{1,p}(M)$, we have

\begin{equation} \label{Riemann-non-sharp-no-tau}
\left( \int_M |u|^r\; dv_g \right)^{\frac{p}{r \theta}} \leq
\left( C \int_M |\nabla_g u|^p\; dv_g + D \int_M |u|^p\; dv_g
\right) \left( \int_M |u|^q\; dv_g \right)^{\frac{p(1 -
\theta)}{\theta q}}\, ,
\end{equation}

\n where $1 < p < n$, $1 \leq q < r < p^*$ and
$\theta = \frac{np(r - q)}{r(q(p - n) + np)} \in (0,1)$ is the
interpolation parameter.

It is a simple matter to add an additional parameter $\tau\geq 1$ in this inequality. We will work with

\begin{gather}\label{AB1}
\left( \int_M |u|^r\; dv_g \right)^{\frac{\tau}{r \theta}} \leq
\left( A \left(\int_M |\nabla_g u|^p\;
dv_g\right)^{\frac{\tau}{p}} + B \left(\int_M |u|^p\;
dv_g\right)^{\frac{\tau}{p}} \right) \left( \int_M |u|^q\; dv_g
\right)^{\frac{\tau(1 - \theta)}{\theta q}}\, , \tag{$GN_R(A,B)$}
\end{gather}
\n  
with $p,q,r,\theta$ as above. Note that when $\tau = p$
we recover (\ref{Riemann-non-sharp-no-tau}). 

The presence of this parameter represents the study of the Gagliardo-Niremberg inequalities 
for the family of equivalent norms in $W^{1,p}$ given by
\[
\left( \left(\int_M |\nabla_g u|^p\;
dv_g\right)^{\frac{\tau}{p}} +  \left(\int_M |u|^p\;
dv_g\right)^{\frac{\tau}{p}} \right)^{1/\tau}
\]
A similar parameter was considered by Druet (\cite{D2})  in the context of Sobolev inequalities 
(his $\theta$ corresponds to our $\tau$), 
in the process of solving  a conjecture of Aubin (conjecture 2 of \cite{Au1}).

 Observe that the non-sharp inequality $GN_R(A,B)$, implies that

\begin{equation}\label{desin}
A \geq A(p,q,r,n)^{\frac{\tau}{p}} \; ,
\end{equation}

\n for any $1 < p < n$ and $1 \leq q < r < p^*$. This is shown by taking an appropriate localized test
function, with support contained in a small enough normal neighbourhood so that the metric is
almost Euclidean, compare \cite{DHV}.

We now study the optimal inequality. Having two constants, the optimality can be defined in two ways. We
follow the more interesting one from the PDE viewpoint (see chapters 4 and 5 \cite{He}):
define the {\em first Riemannian $L^p$-Gagliardo-Nirenberg best constant}
by
\[
A_{opt} = \inf \{ A \in \R:\; \mbox{there exists} \hspace{0,18cm}
B \in \R \hspace{0,18cm} \mbox{such that} \hspace{0,18 cm}
GN_R(A,B) \hspace{0,18cm} \mbox{is valid}\}\, .
\]

\n This optimal constant is positive by (\ref{desin}). Moreover,

\begin{equation}\label{desinf}
A_{opt}^{\frac{p}{\tau}} \geq A(p,q,r,n) \; ,
\end{equation}

\n for any $1 < p < n$ and $1 \leq q < r < p^*$. Then the {\em first optimal Riemannain
$L^p$-Gagliardo-Nirenberg inequality} means that there exists
a constant $B \in \R$ such that, for any $u \in H^{1,p}(M)$,

\[
\left( \int_M |u|^r\; dv_g \right)^{\frac{\tau}{r \theta}} \leq
\left( A_{opt} \left(\int_M |\nabla_g u|^p\;
dv_g\right)^{\frac{\tau}{p}} + B \left(\int_M |u|^p\;
dv_g\right)^{\frac{\tau}{p}} \right) \left( \int_M |u|^q\; dv_g
\right)^{\frac{\tau(1 - \theta)}{\theta q}}\, ,
\]

\n is valid. In contrast with the Euclidean case,
the validity of the optimal inequality is delicate since as $A \to A_{opt}$ the corresponding $B$ might in principle go to infinity. In fact, when $\tau = p >2$  there exists cases where the optimal inequality is not valid,  depending on the geometry of $(M,g)$, (see \cite {CMMZ} and \cite{D5}).

Assuming that $GN_R(A_{opt},B)$ holds, we can define the
{\em second Riemannian $L^p$-Gagliardo-Nirenberg best constant} by

\[
B_{opt} = \inf\{B \in \R; GN_R(A_{opt},B) \hspace{0,2cm} is
\hspace{0,2cm} valid\} \; .
\]

\n Since constant non zero functions belong to $H^{1,p}(M)$,  the constant $B_{opt}$ satisfies

\begin{equation}\label{sco}
B_{opt} \geq |M|^{- \frac{\tau}{n}} \; ,
\end{equation}

\n where  $|M|$ denotes the volume of $(M,g)$.

Then, by the {\em optimal Riemannain $L^p$-Gagliardo-Nirenberg
inequality} we mean that for all $u \in H^{1,p}(M)$,

\[
\left( \int_M |u|^r\; dv_g \right)^{\frac{\tau}{r \theta}} \leq
\left( A_{opt} \left(\int_M |\nabla_g u|^p\;
dv_g\right)^{\frac{\tau}{p}} + B_{opt} \left(\int_M |u|^p\;
dv_g\right)^{\frac{\tau}{p}} \right) \left( \int_M |u|^q\; dv_g
\right)^{\frac{\tau(1 - \theta)}{\theta q}}\, ,
\]

\n is valid. A non-zero function satisfying  equality in $GN_R(A_{opt},B_{opt})$ is called an {\em extremal function}.

We now state the main results of this paper:

\begin{mainthm}\label{tgno1}
Let $(M,g)$ be a smooth, compact Riemannian manifold without
boundary of dimension $n \geq 2$, $1 \leq \tau \leq p$ and $1 \leq
q < r \leq p < n$. If $1 \leq \tau \leq \min\{p,2\}$ then
$GN_R(A_{opt},B)$ is always valid for some $B$.
\end{mainthm}

Theorem \ref{tgno1} complements the results of Chen and Sun \cite{CS}; for $\tau =1$, they deal with the case $p<r$ whereas we study the
case $p\geq r$. The case $p=r$, Nash inequalities, is of particular interest since starting form Nash inequalities one can obtain entropy
inequalities in the same spirit as in \cite{ceccon-montenegro-logarithmic}; the classical case $\tau = p = r = 2$, $q = 1$ is treated by Humbert \cite{Hu}; the conjecture of Aubin 
(conjeture 2 in \cite{Au1}) is set up with  $\tau =\frac{p}{p-1}$. The condition $1 \leq \tau \leq \min\{p,2\}$ extends a similar condition for Sobolev inequalities 
present in Druet's solution \cite{D2} of Aubin's conjecture.

We remark that the arguments used in the proof are of {\em non-local} nature, in contrast with the local techniques used in \cite{CS}; the local arguments being inadequate when $r < p$. These non-local ideas allow the proof of an $L^r$-concentration result (section \ref{concentration}) and a more refined pointwise estimate of certain maximizers (section \ref{pw}). This refinement is essential for the case $r<p$.

Having theorem \ref{tgno1} allows the consideration of the second optimal constant; now by definition  the optimal inequality
 $GN_R(A_{opt},B_{opt})$ holds. We have

\begin{mainthm}\label{extremal}
Let $(M,g)$ be a smooth, compact Riemannian manifold without boundary of
dimension $n \geq 2$, $1 \leq \tau \leq p$ and $1 \leq q < r \leq p < n$. If $1 \leq \tau < 2$ then
$GN_R(A_{opt},B_{opt})$ admits an extremal function.
\end{mainthm}

\section{Proof of Theorem \ref{tgno1}}

We proceed by contradiction: assume inequality $GN_R(A_{opt},B)$ is false for all $B$;
this means that for any
$\alpha>0$ there exists $u \in H^{1,p}(M)$ such that

\begin{equation} \label{contradicao1}
\left( \int_M |u|^r\; dv_g \right)^{\frac{\tau}{r \theta}}\left(
\int_M |u|^q\; dv_g \right)^{-\frac{\tau(1 - \theta)}{\theta q}} -
\alpha \left(\int_M |u|^p\; dv_g\right)^{\frac{\tau}{p}} >
 A_{opt} \left(\int_M |\nabla_g u|^p\;
dv_g\right)^{\frac{\tau}{p}} \, .
\end{equation}

Consider the space $E=\{u \in H^{1,p}(M):\; ||\nabla_g u||_{L^p(M)} = 1\}$. By suitably normalizing, we can assume the
an $u$ satisfying the previous inequality belongs to $E$, that is

\begin{equation} \label{contradicao2}
\left( \int_M |u|^r\; dv_g \right)^{\frac{\tau}{r \theta}}\left(
\int_M |u|^q\; dv_g \right)^{-\frac{\tau(1 - \theta)}{\theta q}} -
\alpha \left(\int_M |u|^p\; dv_g\right)^{\frac{\tau}{p}} >
 A_{opt}  \, .
\end{equation}

Consider now the functional $J_\alpha:E\to \R$ given by the left-hand side of the preceeding inequality. We will
show:

\begin{enumerate}
 \item $J_\alpha$ always admit a maximizer $\tilde u_\alpha$, that is,
$J_\alpha(\tilde u_\alpha) = \sup_{u \in E} J_\alpha =: \nu_\alpha
> A_{opt}$, which satisfies an elliptic PDE as Euler-Lagrange
equation. It will be simpler to work with the normalization
$u_\alpha = \tilde u_\alpha/||\tilde{u}_\alpha||_{L^r(M)}$, which,
by homogeneity, satisfies inequality (\ref{contradicao1}). This
step is studied in section \ref{EL} 

\item Aided by the
Euler-Lagrange equation, satisfied by  $u_\alpha$, we show that it
is concentrated around its maximum  $x_0$ in a sense to be made
precise in section \ref{concentration}, and it also satisfies a
pointwise estimate that quantifies the rate of decay of $u_\alpha$
in terms of the distance to $x_0$; this is done in section
\ref{pw}. 

\item The previous item allows to localize the
integrations in a small normal coordinate chart, and control the non-Euclidean 
terms in the Cartan expansion of the metric-volume-gradient etc. Then using the Euclidean
Gagliardo-Nirenberg inequality will furnish the desired
contradiction.
\end{enumerate}

\subsection{Maximizers and their Euler-Lagrange equations} \label{EL}


For each $\alpha > 0$, consider the functional $J_\alpha$

\[
J_\alpha(u) = \left( \int_M |u|^r\; dv_g\right)^{\frac{\tau}{r
\theta}}\left(\int_M |u|^q dv_g \right)^{- \frac{\tau(1 -
\theta)}{\theta q}} - \alpha \left(\int_M |u|^p\; dv_g
\right)^{\frac{\tau}{p}}
\]

\n defined on the space $E$, and

\begin{equation}\label{3dha}
\nu_\alpha = \sup_{u \in E} J_\alpha (u) > A_{opt}\, .
\end{equation}

\n Note that $\nu_\alpha$ is well-defined and finite  since there are constants $A,B$ such that
$GN_R(A,B)$ holds.

\n Assume first $q > 1$. Since $J_\alpha$ is of class $C^1$, by
using standard variational arguments, we find a maximizer
$\tilde{u}_\alpha \in E$ of $J_\alpha$, i.e.

\begin{equation}\label{3iha}
J_\alpha(\tilde{u}_\alpha) = \nu_\alpha = \sup_{u \in E}
J_\alpha(u)\ .
\end{equation}

\n When $q = 1$, the functional $J_\alpha$ is not $C^1$; however in this case, following Humbert \cite{Hu} we
obtain the existence of extremal satisfying the corresponding Euler-Lagrange equation.
From now on, the arguments are similar in the two cases $q
> 1$ and $q = 1$. Thereby, we will focus our attention only on the case $q > 1$.

\n By (\ref{3iha}), $\tilde{u}_\alpha$ satisfies the Euler-Lagrange
equation

\begin{equation}\label{eq1}
\frac{1}{\theta}||\tilde{u}_\alpha||_{L^r(M)}^{\frac{\tau - r
\theta}{\theta}} ||\tilde{u}_\alpha||_{L^q(M)}^{- \frac{\tau(1 -
\theta)}{\theta}} \tilde{u}_\alpha^{r - 1} - \frac{(1 -
\theta)}{\theta}
||\tilde{u}_\alpha||_{L^r(M)}^{\frac{\tau}{\theta}}
||\tilde{u}_\alpha||_{L^q(M)}^{-\frac{\tau(1 - \theta) + \theta
q}{\theta}} \tilde{u}_\alpha^{q - 1} - \alpha
||\tilde{u}_\alpha||_{L^p(M)}^{\tau - p} \tilde{u}_\alpha^{p - 1}
= \nu_\alpha \Delta_{p,g} \tilde{u}_\alpha \; ,
\end{equation}

\n where $\Delta_{p,g} = -{\rm div}_g(|\nabla_g|^{p-2} \nabla_g)$
is the $p$-Laplace operator of $g$. Because $\nabla_g
|\tilde{u}_\alpha| = \pm \nabla_g \tilde{u}_\alpha$, we can assume
$\tilde{u}_\alpha \geq 0$.

We now set $u_\alpha =
\frac{\tilde{u}_\alpha}{||\tilde{u}_\alpha||_{L^r(M)}}$. Writing the Euler-Lagrange equation in terms of $u_\alpha$, we have

\begin{equation} \label{3ep}
\lambda_\alpha^{-1} A_\alpha \Delta_{p,g} u_\alpha + \alpha
 A_\alpha^{\frac{\tau}{p}} ||u_\alpha||_{L^p(M)}^{\tau -
p} u_\alpha^{p - 1} + \frac{1 - \theta}{\theta}
||u_\alpha||_{L^q(M)}^{-q} u_\alpha^{q - 1} = \frac{1}{\theta}
u_\alpha^{r - 1}\ \ {\rm on}\ \ M\, ,
\end{equation}

\n where  $||u_\alpha||_{L^r(M)} = 1$,

\[
A_\alpha = \left(\int_M u_\alpha^q\; dv_g \right)^{\frac{p(1 -
\theta)}{\theta q}}
\]

\n and

\[
\lambda_\alpha = \nu_\alpha^{-1} ||u_\alpha||_{L^q(M)}^{\frac{(p -
\tau)(1 - \theta)}{\theta}} ||\tilde{u}_\alpha||_{L^r(M)}^{\tau -
p} \; .
\]

\n By Tolksdorf's regularity theory (see \cite{To}), it follows that
$u_\alpha$ is of class $C^1$.

\n We now highlight two important consequences of the Euler-Lagrange
equation.
\begin{equation}\label{rm}
\lambda_\alpha^{-1} \geq A_{opt}^{\frac{p}{\tau}} \; ,
\end{equation}
\n and
 \begin{equation} \label{3}
\lim _{\alpha \rightarrow \infty} A_\alpha = 0\, .
\end{equation}

 In order to show (\ref{rm}), note that taking $\tilde{u}_\alpha$ as test function in (\ref{eq1}), we have
\[
\nu_\alpha \leq
||\tilde{u}_\alpha||_{L^r(M)}^{\frac{\tau}{\theta}}
||\tilde{u}_\alpha||_{L^q(M)}^{- \frac{\tau(1 - \theta)}{\theta}}
\; .
\]

\n Putting together the previous inequality, (\ref{3dha}) and noting that $\tau \leq p$, we get

\[
\lambda_\alpha = \nu_\alpha^{-1}
||\tilde{u}_\alpha||_{L^r(M)}^{\frac{(\tau - p)(1 -
\theta)}{\theta}} ||\tilde{u}_\alpha||_{L^q(M)}^{-\frac{(\tau -
p)(1 - \theta)}{\theta}} ||\tilde{u_\alpha}||_{L^r(M)}^{\tau - p}
= \nu_\alpha^{-1}
\left(||\tilde{u}_\alpha||_{L^r(M)}^{\frac{\tau}{\theta}}
||\tilde{u}_\alpha||_{L^q(M)}^{- \frac{\tau(1 - \theta)}{\theta}}
\right)^{\frac{\tau - p}{\tau}} \leq \nu_\alpha^{-\frac{p}{\tau}}
\leq A_{opt}^{- \frac{p}{\tau}} \; .
\]

The limit (\ref{3}) is shown as follows: first,
\[
1 = \int_M u_\alpha^r dv_g \leq |M|^{1 - \frac{r}{p}} \left(\int_M
u_\alpha^p dv_g \right)^{\frac{r}{p}} \; ,
\]

\n and thus we have $||u_\alpha||_{L^p(M)} > c > 0$ for all $\alpha$. Using
now (\ref{3ep}), we obtain

\begin{equation}\label{dbd}
\alpha A_\alpha^{\frac{\tau}{p}} \left(\int_M u_\alpha^p dv_g
\right)^{\frac{\tau}{p}} \leq c \; ,
\end{equation}

\n which proves claim (\ref{3}).

\subsection{$L^r$-concentration} \label{concentration}

Let $x_\alpha \in M$ be a maximum point of $u_\alpha$, that is,

\begin{equation}\label{3ix}
u_\alpha(x_\alpha) = ||u_\alpha||_{L^\infty(M)}\, .
\end{equation}

\n Throughout this section, we use the notation
$\lim_{\sigma,\alpha \rightarrow \infty}$ to mean $\lim_{\sigma
\rightarrow \infty} \lim_{\alpha \rightarrow \infty}$.

\n Our aim here is to establish that

\begin{equation} \label{4}
\lim_{\sigma,\alpha \rightarrow \infty} \int_{B(x_\alpha,\sigma
a_\alpha)} u_\alpha^r\; dv_g = 1\, ,
\end{equation}

\n where
\begin{equation}\label{ic}
a_\alpha = A_\alpha^{\frac{r}{np - nr + pr}}\, .
\end{equation}

Let us proceed. Since $np - nr + pr > 0$ and by (\ref{3}), $a_\alpha
\rightarrow 0$ as $\alpha \rightarrow \infty$.

In order to work in local coordinates, we pullback the metric and the function $u_\alpha$ to a ball $B(0,\sigma)$ in 
$T_{x_\alpha}M$. The pullbacks will be by exponential map precomposed by a dilation: if $\delta_L:TM \to TM$ denotes the dilation by $L$, $\delta_L(x) = Lx$ and $E_\alpha=\exp_{x_\alpha}\circ \delta_{a_\alpha}$, define 
\begin{equation}\label{l}
\begin{array}{l}
h_\alpha = E_\alpha^*g\, , \vspace{0,3cm}\\
\varphi_\alpha = a_\alpha^{\frac{n}{r}}
u_\alpha\circ E_\alpha  \, .
\end{array}
\end{equation}

Note that since $a_\alpha \rightarrow 0$ as $\alpha \rightarrow \infty$, these are well defined: for $\alpha$ large 
enough, $B(0,a_\alpha\sigma)$ will be contained inside of a set where $\exp_{x_\alpha}$ is a diffeomorphism.

\n By (\ref{3ep}), one easily deduces that

\begin{equation}\label{eqgn}
\lambda_\alpha^{-1} \Delta_{p,h_\alpha} \varphi_\alpha + \alpha
A_\alpha^{\frac{\tau}{p} - 1} ||u_\alpha||_{L^p(M)}^{\tau - p}
a_\alpha^p \varphi_\alpha^{p - 1} + \frac{1 - \theta}{\theta}
\varphi_\alpha^{q - 1} = \frac{1}{\theta} \varphi_\alpha^{r - 1}\
\ {\rm on}\ \ B(0, \sigma)\, .
\end{equation}

\n By (\ref{rm}) and applying the Moser's iterative scheme (see
\cite{Mo}) to this last equation, we see that, for $\alpha$ large enough,

\[
a_\alpha^{n} ||u_\alpha||_{L^\infty(M)}^r = \sup_{B(0,\frac{\sigma}{2})} \varphi_\alpha^r \leq c \int_{B(0,\sigma)} \varphi_\alpha^r\; dh_\alpha =
c \int_{B(x_\alpha, \sigma a_\alpha)} u_\alpha^r\; dv_g \leq c
\]

\n  This estimate together with
\[
1 = \int_M u_\alpha^r\; dv_g \leq ||u_\alpha||_{L^\infty(M)}^{r-q} \int_M u_\alpha^q dv_g = \left(||u_\alpha||_{L^\infty(M)} \;
a_\alpha^{\frac{n}{r}}\right)^{r-q}
\]

\n gives

\begin{equation}\label{3dqp}
1 \leq ||u_\alpha||_{L^\infty(M)} \; a_\alpha^{\frac{n}{r}} \leq c \, .
\end{equation}

\n We will use in the sequel that (\ref{3dqp}) means that the
limiting behaviour of $||u_\alpha||_{L^\infty(M)} $ and $
a_\alpha^{-\frac{n}{r}}$ have the same order.

 In particular, there exists a constant $c > 0$ such that

\begin{equation}\label{lb}
\int_{B(0,\sigma)} \varphi_\alpha^r\; dh_\alpha \geq c > 0
\end{equation}

\n for $\alpha$ large enough.

\n Now using Cartan expansion in normal coordinates and (\ref{3dqp}), we have for each $\sigma > 0$, that

\[
\int_{B(0,\sigma)} \varphi_\alpha^p dx \leq c \int_{B(0,\sigma)}
\varphi_\alpha^p dh_\alpha = c \; a_\alpha^{\frac{n}{r}(p - r)}
\int_{B(x_\alpha, \sigma a_\alpha)} u_\alpha^p dv_g \leq c
a_\alpha^{\frac{n}{r}(p - r)}||u_\alpha||^p_{L^\infty(M)} \sigma^n
a_\alpha^n \leq c(\sigma) \; ,
\]

\n where $c(\sigma) \rightarrow \infty$ when $\sigma \rightarrow
\infty$, but $c(\sigma)$ being independent of $\alpha$. We have also

\[
\int_{B(0,\sigma)} |\nabla \varphi_\alpha|^p\; dx \leq c
\int_{B(0,\sigma)} |\nabla_{h_\alpha} \varphi_\alpha|^p\;
dh_\alpha = c A_\alpha \int_{B(x_\alpha,\sigma a_\alpha)}
|\nabla_g u_\alpha|^p\; dv_g  \leq c \;
A_{opt}^{-\frac{p}{\tau}}\, .
\]

\n Therefore there exists $\varphi \in W^{1,p}(\R^n)$ such that, for some
subsequence, $\varphi_\alpha \rightharpoonup \varphi$ in
$W^{1,p}_{loc}(\R^n)$. For each $\sigma > 0$, we have

\[
\int_{B(0,\sigma)} \varphi^q\; dx = \lim_{\alpha \rightarrow
\infty} \int_{B(0,\sigma)} \varphi_\alpha^q\; dh_\alpha =
\lim_{\alpha \rightarrow \infty} \frac{\int_{B(x_\alpha,\sigma
a_\alpha)} u_\alpha^q\; dv_g}{\int_M u_\alpha^q\; dv_g} \leq 1
\]

\n and

\[
\int_{B(0,\sigma)} \varphi^r\; dx = \lim_{\alpha \rightarrow
\infty} \int_{B(0,\sigma)} \varphi_\alpha^r\; dh_\alpha =
\lim_{\alpha \rightarrow \infty} \int_{B(x_\alpha,\sigma
a_\alpha)} u_\alpha^r\; dv_g \leq 1 \; .
\]

\n In particular,

\begin{equation}\label{inter}
\varphi \in L^q(\R^n) \cap L^r(\R^n) \; .
\end{equation}

\n Let $\eta \in C_0^1(\R)$ be a cutoff function such that $\eta =
1$ on $[0,\frac{1}{2}]$, $\eta = 0$ on $[1,\infty)$ and $0 \leq
\eta \leq 1$. Define $\eta_{\alpha,\sigma}(x) = \eta((\sigma
a_\alpha)^{-1} d_g(x,x_\alpha))$. Choosing $u_\alpha
\eta_{\alpha,\sigma}^p$ as a test function in (\ref{3ep}), one
gets

\[
\lambda_\alpha^{-1} A_\alpha \int_M |\nabla_g u_\alpha|^p
\eta_{\alpha,\sigma}^p\; dv_g + \lambda_\alpha^{-1}A_\alpha \int_M
|\nabla_g u_\alpha|^{p - 2} \nabla_g u_\alpha \cdot \nabla_g
(\eta_{\alpha,\sigma}^p) u_\alpha\; dv_g + \frac{1 -
\theta}{\theta}  \frac{\int_M u_\alpha^q \eta_{\alpha,\sigma}^p\;
dv_g}{\int_M u_\alpha^q dv_g}
\]

\begin{equation}\label{3sc3}
\leq \frac{1}{\theta} \int_M u_\alpha^r \eta_{\alpha,\sigma}^p\;
dv_g\, .
\end{equation}

\n We now show that

\begin{equation}\label{ef1}
\lim_{\sigma,\alpha \rightarrow \infty} A_\alpha \int_M |\nabla_g
u_\alpha|^{p - 2} \nabla_g u_\alpha \cdot \nabla_g
(\eta_{\alpha,\sigma}^p) u_\alpha\; dv_g = 0\, .
\end{equation}

\n Taking $u_\alpha$ as test function, by (\ref{3ep}) we have

\[
A_\alpha \int_M |\nabla_gu_\alpha|^p dv_g \leq \lambda_\alpha \leq
A_{opt}^{-\frac{p}{\tau}} \; .
\]

\n Therefore, it suffices to establish that

\begin{equation}\label{ef6}
A_\alpha \int_M u_\alpha^p |\nabla_g \eta_{\alpha,\sigma}|^p\;
dv_g \leq \frac{c}{\sigma^p}\, .
\end{equation}

\n Using  (\ref{ic}) and (\ref{3dqp}), we derive

\[
A_\alpha \int_M u_\alpha^p |\nabla_g \eta_{\alpha,\sigma}|^p\;
dv_g \leq c \frac{A_\alpha ||u_\alpha||_{L^\infty(M)}^{p -
r}}{\sigma^p a_\alpha^p} \int_M u_\alpha^r\; dv_g \leq c
\frac{A_\alpha a_\alpha^{-\frac{n}{r}(p - r)}}{\sigma ^p
a_\alpha^p} \leq \frac{c}{\sigma^p}.
\]

\n Therefore (\ref{ef6}) holds and (\ref{ef1}) is valid.

\n Replacing (\ref{rm}) and (\ref{ef1}) in (\ref{3sc3}) and
because $\lambda_\alpha^{-1}$ is limited, one arrives at

\[
\theta \lim_{\sigma,\alpha \rightarrow
\infty}\left(A_{opt}^{\frac{p}{\tau}}\; A_\alpha \int_M |\nabla_g
u_\alpha|^p \eta_{\alpha, \sigma}^p\; dv_g \right) + (1 - \theta)
\lim_{\sigma, \alpha \rightarrow \infty} \frac{\int_M u_\alpha^q
\eta_{\alpha,\sigma}^p\; dv_g}{\int_M u_\alpha^q\; dv_g} \leq
\lim_{\sigma, \alpha \rightarrow \infty} \int_M u_\alpha^r
\eta_{\alpha, \sigma}^p\; dv_g\, .
\]

\n To rewrite this inequality in a suitable format, we first remark that

\[
\left| \frac{\int_M u_\alpha^q \eta_{\alpha,\sigma}^p\;
dv_g}{\int_M u_\alpha^q\; dv_g} - \frac{\int_M u_\alpha^q
\eta_{\alpha,\sigma}^q\; dv_g}{\int_M u_\alpha^q\; dv_g} \right|
\leq \frac{\int_{B(x_\alpha,\sigma a_\alpha) \setminus B(x_\alpha,
\sigma a_\alpha/2)} u_\alpha^q\; dv_g}{\int_M u_\alpha^q\; dv_g} =
\int_{B(0,\sigma) \setminus B(0,\sigma/2)} \varphi_\alpha^q\;
dh_\alpha
\]

\n and the above right-hand side converges to $0$ as $\sigma
\rightarrow \infty$, by (\ref{inter}). So,

\[
\lim_{\sigma,\alpha \rightarrow \infty} \frac{\int_M u_\alpha^q
\eta_{\alpha,\sigma}^p\; dv_g}{\int_M u_\alpha^q\; dv_g} =
\lim_{\sigma,\alpha \rightarrow \infty} \frac{\int_M u_\alpha^q
\eta_{\alpha,\sigma}^q\; dv_g}{\int_M u_\alpha^q\; dv_g}\, .
\]

\n In a similar way,

\[
\left| \int_M u_\alpha^r \eta_{\alpha,\sigma}^p \; dv_g - \int_M
u_\alpha^r \eta_{\alpha,\sigma}^r \;dv_g \right| \leq
\int_{B(x_\alpha,\sigma a_\alpha) \setminus B(x_\alpha, (\sigma
a_\alpha)/2)} u_\alpha^r \;dv_g = \int_{B(0,\sigma) \setminus
B(0,\sigma/2)} \varphi_\alpha^r dh_\alpha \; ,
\]

\n using (\ref{inter}), so that

\[
\lim_{\sigma,\alpha \rightarrow \infty} \int_M u_\alpha^r
\eta_{\alpha,\sigma}^r\; dv_g = \lim_{\sigma,\alpha \rightarrow
\infty} \int_M u_\alpha^r \eta_{\alpha,\sigma}^p\; dv_g\, .
\]

\n Consequently, we can write

\begin{equation}\label{ef2}
\theta \lim_{\sigma,\alpha \rightarrow
\infty}\left(A_{opt}^{\frac{p}{\tau}} \; A_\alpha \int_M |\nabla_g
u_\alpha|^p \eta_{\alpha, \sigma}^p\; dv_g \right) + (1 - \theta)
\lim_{\sigma, \alpha \rightarrow \infty} \frac{\int_M u_\alpha^q
\eta_{\alpha,\sigma}^q\; dv_g}{\int_M u_\alpha^q\; dv_g} \leq
\lim_{\sigma, \alpha \rightarrow \infty} \int_M u_\alpha^r
\eta_{\alpha, \sigma}^r\; dv_g\, .
\end{equation}

\n On the other hand, for $\varepsilon > 0$ let the constant
$B_\varepsilon > 0$, independent of $\alpha$, such that

\[
\left( \int_M u_\alpha^r \eta^r_{\alpha,\sigma}\; dv_g
\right)^{\frac{\tau}{r \theta}} \leq \left( (A_{opt} +
\varepsilon) \left(\int_M |\nabla_g(u_\alpha
\eta_{\alpha,\sigma})|^p\; dv_g\right)^{\frac{\tau}{p}} +
B_\varepsilon \left(\int_M u_\alpha^p \eta_{\alpha,\sigma}^p\;
dv_g\right)^{\frac{\tau}{p}} \right) \left(\int_M u_\alpha^q
\eta^q_{\alpha,\sigma}\; dv_g \right)^{\frac{\tau(1 -
\theta)}{\theta q}}\, .
\]

\n From the definition of $A_\alpha$, Young inequality and $(x + y)^p \leq x^p + c x^{p - 1}y + c y^p$ for $x,y \geq 0$, one
has

\[
\left(\int_M u_\alpha^r \eta^r_{\alpha,\sigma} dv_g
\right)^{\frac{\tau}{r \theta}} \leq c \left(A_\alpha \int_M
u_\alpha^p dv_g \right)^{\frac{\tau}{p}}
\]
\[
+ (A_{opt} + \varepsilon) \left(\left((1 + \varepsilon) \int_M
|\nabla_g u_\alpha|^p \eta_{\alpha,\sigma}^p\; dv_g +
c(\varepsilon) \int_M u_\alpha^p |\nabla_g \eta_{\alpha,\sigma}|^p
dv_g \right)^{\frac{\tau}{p}}\right) \left(\int_M u_\alpha^q
\eta_{\alpha,\sigma}^q\; dv_g \right)^{\frac{\tau(1 -
\theta)}{\theta q}} \, .
\]

\n Then, using (\ref{dbd}), (\ref{ef6}), letting $\alpha,\sigma
\rightarrow \infty$ and $\varepsilon \rightarrow 0$, one gets

\begin{equation}\label{3il2}
\lim_{\sigma,\alpha \rightarrow \infty}\left( \int_M u_\alpha^r
\eta_{\alpha,\sigma}^r\; dv_g \right)^{\frac{p}{r \theta}} \leq
\lim_{\sigma,\alpha \rightarrow \infty} \left(
A_{opt}^{\frac{p}{\tau}} A_\alpha \int_M |\nabla_g u_\alpha|^p
\eta_{\alpha,\sigma}^p\; dv_g \right) \lim_{\sigma,\alpha
\rightarrow \infty} \left(\frac{\int_M u_\alpha^q
\eta_{\alpha,\sigma}^q\; dv_g}{\int_M u_\alpha^q\; dv_g}
\right)^{\frac{p(1 - \theta)}{\theta q}}\, .
\end{equation}

\n Let

\[
X = \lim_{\sigma, \alpha \rightarrow \infty} \left(
A_{opt}^{\frac{p}{\tau}} \; A_\alpha \; \int_M |\nabla_g
u_\alpha|^p \eta_{\alpha,\sigma}^p dv_g \right), \ \ Y =
\lim_{\sigma,\alpha \rightarrow \infty}\frac{\int_M u_\alpha^q
\eta_{\alpha,\sigma}^q dv_g}{\int_M u_\alpha^q dv_g}, \ \ Z =
\lim_{\sigma,\alpha \rightarrow \infty} \int_M u_\alpha^r
\eta_{\alpha,\sigma}^r dv_g\, .
\]

\n Clearly, $X, Y, Z \leq 1$ and (\ref{ef2}) and (\ref{3il2}) may be rewritten as

\begin{equation}\label{3sc4}
\left\{%
\begin{array}{c}
\theta X + (1 - \theta) Y \leq Z   \vspace{0,5cm}\\
Z \leq X^{\frac{r \theta}{p}} Y^{\frac{r(1 - \theta)}{q}} \, .
\end{array}%
\right.
\end{equation}

\n By (\ref{lb}), one also has $Z > 0$, so that $X, Y > 0$.

\n Assertion (\ref{4}) follows readily by proving that $Z =
1$. For this, we will consider the behavior of $u_\alpha$ outside
$B(x_\alpha, \sigma a_\alpha)$. Indeed, let $\zeta_{\alpha,\sigma}
= 1 - \eta_{\alpha,\sigma}$ on $M$. Inequalities (\ref{ef2}) and (\ref{3il2}) remain valid for
$\zeta_{\alpha,\sigma}$ in place of $\eta_{\alpha,\sigma}$. In
other words,

\[
\theta \lim_{\sigma,\alpha \rightarrow \infty}
\left(A_{opt}^{\frac{p}{\tau}} \;A_\alpha \int_M |\nabla_g
u_\alpha|^p \zeta_{\alpha, \sigma}^p\; dv_g\right)  + (1 - \theta)
\lim_{\sigma, \alpha \rightarrow \infty} \frac{\int_M u_\alpha^q
\zeta_{\alpha,\sigma}^q\; dv_g}{\int_M u_\alpha^q\; dv_g} \leq
\lim_{\sigma, \alpha \rightarrow \infty} \int_M u_\alpha^r
\zeta_{\alpha, \sigma}^r\; dv_g
\]

\n and

\[
\lim_{\sigma,\alpha \rightarrow \infty}\left( \int_M u_\alpha^r
\zeta_{\alpha,\sigma}^r\; dv_g \right)^{\frac{p}{r \theta}} \leq
\lim_{\sigma,\alpha \rightarrow \infty} \left(
A_{opt}^{\frac{p}{\tau}} A_\alpha \int_M |\nabla_g u_\alpha|^p
\zeta_{\alpha,\sigma}^p\; dv_g \right) \lim_{\sigma,\alpha
\rightarrow \infty} \left(\frac{\int_M u_\alpha^q
\zeta_{\alpha,\sigma}^q\; dv_g}{\int_M u_\alpha^q\; dv_g}
\right)^{\frac{p(1 - \theta)}{\theta q}}\, .
\]

\n In a similar way, we denote

\[
\tilde{X} = \lim_{\sigma,\alpha \rightarrow \infty}
\left(A_{opt}^{\frac{p}{\tau}} \; A_\alpha \int_M |\nabla_g
u_\alpha|^p \zeta_{\alpha, \sigma}^p\; dv_g\right) , \ \ \tilde{Y}
= \lim_{\sigma,\alpha \rightarrow \infty}\frac{\int_M u_\alpha^q
\zeta_{\alpha,\sigma}^q\; dv_g}{\int_M u_\alpha^q\; dv_g},\ \
\tilde{Z} = \lim_{\sigma,\alpha \rightarrow \infty} \int_M
\zeta_{\alpha,\sigma}^r u_\alpha^r\, dv_g\, ,
\]

\n so that

\begin{equation}\label{ef3}
\left\{%
\begin{array}{c}
  \theta \tilde{X} + (1 - \theta) \tilde{Y} \leq \tilde{Z}  \vspace{0,5cm}\\
  \tilde{Z} \leq \tilde{X}^{\frac{r \theta}{p}} \tilde{Y}^{\frac{r(1 - \theta)}{q}}  \, .\\
\end{array}%
\right.
\end{equation}

\n We next assert that

\begin{equation} \label{s1}
Y + \tilde{Y} = 1 \ \ \mbox{and} \ \ Z + \tilde{Z} = 1\, .
\end{equation}

\n To justify the first equality, let us write

\[
1 = \frac{\int_M u_\alpha^q \eta_{\alpha,\sigma}^q\; dv_g}{\int_M
u_\alpha^q\; dv_g} + \frac{\int_M u_\alpha^q (1 -
\eta_{\alpha,\sigma}^q)\; dv_g}{\int_M u_\alpha^q\; dv_g} \, .
\]

\n Then,

\[
\left| \frac{\int_M u_\alpha^q (1 - \eta_{\alpha,\sigma}^q)\; dv_g
- \int_M u_\alpha^q \zeta_{\alpha,\sigma}^q\; dv_g}{\int_M
u_\alpha^q\; dv_g}\right| \leq
\left|\frac{\int_{B(x_\alpha,a_\alpha \sigma) \setminus
B(x_\alpha,a_\alpha \sigma/2)} u_\alpha^q dv_g}{\int_M u_\alpha^q
dv_g} \right| = \int_{B(0,\sigma) \setminus B(0,\sigma/2)}
\varphi_\alpha^q\; dh_\alpha \; .
\]

\n Since $\varphi \in L^q(\R^n)$, we get

\[
\lim_{\sigma,\alpha \rightarrow \infty} \frac{\int_M u_\alpha^q \zeta_{\alpha,\sigma}^q dv_g}{\int_M u_\alpha^q dv_g} = \lim_{\sigma,\alpha \rightarrow \infty} \frac{\int_M u_\alpha^q (1 - \eta_{\alpha,\sigma}^q) dv_g}{\int_M u_\alpha^q dv_g} \; ,
\]

\n which yields $Y + \tilde{Y} = 1$. Analogously, since $\varphi \in
L^r(\R^n)$,

\[
1 = \int_M u_\alpha^r \eta_{\alpha,\sigma}^r\; dv_g + \int_M
u_\alpha^r (1 - \eta_{\alpha,\sigma}^r)\; dv_g
\]

\n and

\[
\left|\int_M u_\alpha^r (1 - \eta_{\alpha,\sigma}^r)\; dv_g -
\int_M u_\alpha^r \zeta_{\alpha,\sigma}^r\; dv_g\right| \leq
\int_{B(x_\alpha,a_\alpha \sigma) \setminus B(x_\alpha,a_\alpha
\sigma/2)} u_\alpha^r\; dv_g = \int_{B(0,\sigma)\setminus
B(0,\sigma/2)} \varphi_\alpha^r\; dh_\alpha
\]

\n which gives $Z + \tilde{Z} = 1$.

\n We  are now ready to prove that $Z = 1$. The first inequality in
(\ref{3sc4}) lead us to three possible alternatives:

\vspace{0,4cm}

\centerline{{(a)} $X \leq Z$ and $Y \leq Z$, \ \ \ \ {(b)} $Y \leq
Z \leq X$, \ \ \ \ {(c)} $X \leq Z \leq Y$\, .}

\vspace{0,4cm}

\n If (a) holds , the second inequality in (\ref{3sc4})
implies

\[
1 \leq Z^{\frac{r \theta}{p} + \frac{r(1 - \theta)}{q} - 1}\, .
\]

\n But the definition of $\theta$ furnishes $\frac{r \theta}{p} +
\frac{r(1 - \theta)}{q} - 1 > 0$, so that $Z = 1$.

\n Suppose then that the item (b) holds. Again, by (\ref{3sc4}),

\[
Z \geq \theta X + (1 - \theta) Y = \frac{r \theta}{p} X + (1 -
\frac{r\theta}{p}) Y + \theta(1 - \frac{r}{p}) X + \theta(-1 +
\frac{r}{p}) Y\, .
\]

\n Using the assumption $p \geq r$, $Y \leq X$ and Young's
inequality, we derive

\[
Z \geq X^{\frac{r \theta}{p}} Y^{\frac{p - r\theta}{p}} + \theta
\frac{p - r}{p} (X - Y) \geq X^{\frac{r \theta}{p}} Y^{\frac{p - r
\theta}{p}}\, ,
\]

\n so that the second inequality in (\ref{3sc4}) immediately
yields

\[
1 \leq Y^{\frac{r(1 - \theta)}{q} - \frac{p - r\theta}{p}}\, .
\]

\n Since $\frac{r(1 - \theta)}{q} - \frac{p - r\theta}{p} > 0$,
one has $Y = 1$. Thus, evoking (\ref{ef3}) and (\ref{s1}), one
easily deduces that $Z = 1$.

\n Finally, we come to the alternative (c). We first show that $Y
= 1$. Otherwise, by (\ref{s1}), one has $\tilde{Y} > 0$ and this
implies, by (\ref{ef3}) and $\theta < 1$, that $\tilde{X} > 0$ and $\tilde{Z} > 0$.
Since, by hypothesis, $Z \leq Y$, again thanks to (\ref{s1}), one
concludes that $\tilde{Y} \leq \tilde{Z}$. Therefore, applying the
previously discussed cases (a) and (b) with $\tilde{X}$,
$\tilde{Y}$ and $\tilde{Z}$ in the place of $X$, $Y$ and $Z$, one
arrives at $\tilde{Z} = 1$. But this contradicts the fact that $Z
> 0$, so that $Y = 1$. As before, (\ref{ef3}) and (\ref{s1})
produce $Z = 1$ and this finishes the proof of the
$L^r$-concentration.

\subsection{Pointwise estimates} \label{pw}

The aim of this section is to prove the following pointwise estimate of the decay of $u_\alpha$ in terms of the
distance to its maximum:

\n \emph{For any constant $\lambda
> 0$ there exists a constant $c_\lambda
> 0$, independent of $\alpha$, such that}

\[
d_g(x,x_\alpha)^\lambda u_\alpha(x) \leq c_\lambda \;
a_\alpha^{\lambda - \frac{n}{r}}
\]

\n \emph{for all $x \in M$ and $\alpha$ large enough.}

\n The proof proceeds by contradiction. Suppose that the assertion above is false.
Then, there exist $\lambda_0 > 0$ and $y_\alpha \in M$ such that
$f_{\alpha}(y_\alpha) \rightarrow \infty$ as $\alpha \rightarrow
\infty$, where

\[
f_{\alpha}(x) = d_g(x,x_\alpha)^{\lambda_0} u_\alpha(x) \;
a_\alpha^{-\lambda_0 + \frac{n}{r}}\, .
\]

\n Assume, without loss of generality, that $f_{\alpha}(y_\alpha) = ||f_{\alpha}||_{L^\infty(M)}$.
From (\ref{3dqp}), we have

\[
f_\alpha(y_\alpha) \leq c \;
\frac{u_\alpha(y_\alpha)}{||u_\alpha||_{L^\infty(M)}}
d_g(x_\alpha,y_\alpha)^{\lambda_0} a_\alpha^{- \lambda_0} \leq c
d_g(x_\alpha,y_\alpha)^{\lambda_0} a_\alpha^{-\lambda_0}\, ,
\]

\n so that
\begin{equation}\label{ldi}
d_g(x_\alpha,y_\alpha) a_\alpha^{-1}
\rightarrow \infty\, .
\end{equation}

\n For any fixed $\sigma > 0$ and $\varepsilon \in (0,1)$, we next show that

\begin{equation}\label{int}
B(y_\alpha,\varepsilon d_g(x_\alpha,y_\alpha)) \cap B(x_\alpha, \sigma a_\alpha) = \emptyset
\end{equation}

\n for $\alpha$ large enough. Clearly, this assertion follows from

\[
d_g(x_\alpha,y_\alpha) \geq \sigma a_\alpha + \varepsilon d(x_\alpha,y_\alpha) \, .
\]

\n But the above inequality is equivalent to
\[
d_g(x_\alpha, y_\alpha)(1 - \varepsilon) a_\alpha^{-1} \geq \sigma \, ,
\]

\n which is clearly satisfied,
since $d_g(x_\alpha,y_\alpha) a_\alpha^{-1} \rightarrow \infty$ and $1 - \varepsilon > 0$.

\n We claim that exists a constant $c > 0$ such that

\begin{equation}\label{lim}
u_\alpha(x) \leq c u_\alpha(y_\alpha)
\end{equation}

\n for all $x \in B(y_\alpha, \varepsilon d_g(x_\alpha,y_\alpha))$
and $\alpha$ large enough. In fact, for each $x \in
B(y_\alpha, \varepsilon d_g(x_\alpha,y_\alpha))$, we have

\[
d_g(x,x_\alpha) \geq d_g(x_\alpha,y_\alpha) - d_g(x,y_\alpha) \geq (1 - \varepsilon)
d_g(x_\alpha,y_\alpha)\, .
\]

\n Thus,

\[
d_g(y_\alpha,x_\alpha)^{\lambda_0} u_\alpha(y_\alpha) a_\alpha
^{-\lambda_0 + \frac{n}{r}} = f_{\alpha}(y_\alpha) \geq
f_{\alpha}(x) = d_g(x,x_\alpha)^{\lambda_0} u_\alpha(x)
a_\alpha^{-\lambda_0 + \frac{n}{r}}
\]

\[
\geq (1 - \varepsilon)^{\lambda_0}
d_g(y_\alpha,x_\alpha)^{\lambda_0} u_\alpha(x)
a_\alpha^{-\lambda_0 + \frac{n}{r}}\, ,
\]

\n so that

\[
u_\alpha(x) \leq \left(\frac{1}{1 -
\varepsilon}\right)^{\lambda_0} u_\alpha(y_\alpha)
\]

\n for all $x \in B(y_\alpha, \varepsilon d_g(x_\alpha,y_\alpha))$ and $\alpha$ large enough. This proves our claim.

\n Since $p \geq r$, by (\ref{ic}) and (\ref{3dqp}), one has

\[
A_\alpha^{\frac{1}{p}} u_\alpha(y_\alpha)^{\frac{p - r}{p}} \leq
A_\alpha^{\frac{1}{p}} ||u_\alpha||_{L^\infty(M)}^{\frac{p -
r}{p}} \leq c \; a_\alpha \rightarrow 0\, .
\]

\n So, in the same way as in (\ref{l}), we can now define, inside of  $B(0,2)$ and for $\alpha$ large enough,
$ E_\alpha = \exp_{y_\alpha}\circ \delta_{A_\alpha^{\frac{1}{p}}
u_\alpha(y_\alpha)^{\frac{p-r}{p}}} $, and pullback
\[
\begin{array}{l}
h_\alpha =  E_\alpha^\star g \\
\psi_\alpha= u_\alpha(y_\alpha)^{-1}
u_\alpha \circ E_\alpha 
\end{array}
\]

\n From (\ref{3ep}), it readily follows that

\begin{equation}\label{3sc5}
\lambda_\alpha^{-1} \Delta_{p,h_\alpha} \psi_\alpha + \alpha
A_\alpha^{\frac{\tau}{p}} ||u_\alpha||_{L^p(M)}^{\tau - p}
u_\alpha(y_\alpha)^{p - r} \psi_\alpha^{p - 1} + \frac{1 -
\theta}{\theta} ||u_\alpha||_{L^q(M)}^{-q} u_\alpha(y_\alpha)^{q -
r} \psi_\alpha^{q - 1} = \frac{1}{\theta} \psi_\alpha^{r - 1} \ \
{\rm on}\ \ B(0,2)\, .
\end{equation}

\n In particular,

\[
\int_{B(0,2)} |\nabla_{h_\alpha} \psi_\alpha|^{p - 2}
\nabla_{h_\alpha} \psi_\alpha \cdot  \nabla_{h_\alpha} \phi\;
dv_{h_\alpha} \leq c \int_{B(0,2)} \psi_\alpha^{r - 1} \phi\;
dv_{h_\alpha}
\]

\n for all positive test function $\phi \in C_0^1(B(0,2))$. So, by the Moser's iterative scheme and (\ref{3dqp}), one deduces that
\[
1 \leq \sup_{B(0,\frac{1}{4})} \psi_\alpha^r \leq c
\int_{B(0,\frac{1}{2})} \psi_\alpha^r\; dv_{h_\alpha} = c \left(
A_\alpha^{\frac{\theta q}{p(1 - \theta)}} u_\alpha(y_\alpha)^{r -
q} \right)^{-\frac{n(1 - \theta)}{\theta q}} \int_{B(y_\alpha,
\frac{1}{2} A_\alpha^{\frac{1}{p}}
u_\alpha(y_\alpha)^{\frac{p-r}{p}})} u_\alpha^r\; dv_g
\]
\[
\leq c \left(\frac{||u_\alpha||_{L^\infty(M)}}{u_\alpha(y_\alpha)}
\right)^{\frac{np - rn + pr}{p}} \int_{B(y_\alpha, \frac{1}{2}
A_\alpha^{\frac{1}{p}} u_\alpha(y_\alpha)^{\frac{p-r}{p}})}
u_\alpha^r\; dv_g\,.
\]

\n For simplicity, we rewrite this last inequality as

\begin{equation}\label{3dcm}
0 < c \leq m_\alpha^\varrho \int_{B(y_\alpha,\frac{1}{2}
A_\alpha^{\frac{1}{p}} u_\alpha(y_\alpha)^{\frac{p-r}{p}})}
u_\alpha^r\; dv_g\, ,
\end{equation}

\n where $m_\alpha =
\frac{||u_\alpha||_{L^\infty(M)}}{u_\alpha(y_\alpha)}$ and
$\varrho = \frac{np - rn + pr}{p} > 0$.

\n By (\ref{ic}), (\ref{3dqp}) and (\ref{ldi}), note that
$B(y_\alpha, \frac{1}{2} A_\alpha^{\frac{1}{p}}
u_\alpha(y_\alpha)^{\frac{p - r}{p}}) \subset
B(y_\alpha,\varepsilon d(x_\alpha,y_\alpha))$ for $\alpha$ large enough. So, the $L^r$-concentration
property (\ref{4}) combined with (\ref{int}) provide

\[
\int_{B(y_\alpha, \frac{1}{2} A_\alpha^{\frac{1}{p}}
u_\alpha(y_\alpha)^{\frac{p-r}{p}})} u_\alpha^r\; dv_g \rightarrow
0 \; ,
\]

\n when $\alpha \rightarrow \infty$.

\n So, we have that

\[
\lim_{\alpha \rightarrow \infty} m_\alpha = \infty \;.
\]

\n Our goal now is to establish a contradiction to (\ref{3dcm}). Initially, from (\ref{3dqp}) and (\ref{lim}), we have

\begin{equation}\label{3dpt}
m_\alpha^\varrho \int_{D_{\alpha}} u_\alpha^r\; dv_g \leq
m_\alpha^\varrho ||u_\alpha||^r_{L^\infty(D_{\alpha})}
(A_\alpha^{\frac{1}{p}} u_\alpha(y_\alpha)^{\frac{p-r}{p}})^n \leq
c m_\alpha^\varrho u_\alpha(y_\alpha)^r (A_\alpha^{\frac{1}{p}}
u_\alpha(y_\alpha)^{\frac{p-r}{p}})^n \leq c\, ,
\end{equation}

\n where $D_{\alpha} = B(y_\alpha, A_\alpha^{\frac{1}{p}}
u_\alpha(y_\alpha)^{\frac{p-r}{p}})$ \; .

\n Consider the function $\eta_\alpha(x) = \eta( A_\alpha^{-\frac{1}{p}} d_g(x,y_\alpha) u_\alpha(y_\alpha)^{\frac{r-p}{p}})$, where $\eta
\in C_0^1(\R)$ is a cutoff function such that $\eta = 1$ on $[0,\frac{1}{2}]$, $\eta = 0 $ on $[1,\infty)$ and $0 \leq \eta
\leq 1$. Taking $u_\alpha \eta_\alpha^p$ as a test function in (\ref{3ep}), one has

\[
\lambda_\alpha ^{-1} A_\alpha \int_M |\nabla_g u_\alpha|^p
\eta_\alpha^p\; dv_g + \lambda_\alpha^{-1} p A_\alpha \int_M |\nabla_g u_\alpha|^{p -
2} u_\alpha \eta_\alpha^{p - 1} \nabla_g u_\alpha \cdot \nabla_g
\eta_\alpha\; dv_g + \alpha A_\alpha^{\frac{\tau}{p}}
||u_\alpha||_{L^p(M)}^{\tau - p} \int_M u_\alpha^p \eta_\alpha^p\;
dv_g
\]
\[
+ \frac{1 - \theta}{\theta} \frac{\int_M u_\alpha^q
\eta_\alpha^p}{\int_M u_\alpha^p dv_g}\; dv_g = \frac{1}{\theta}
\int_M u_\alpha^r \eta_\alpha^p\; dv_g\, .
\]

\n By H\"{o}lder and Young's inequalities,
\[
\left| \int_M |\nabla_g u_\alpha|^{p - 2} u_\alpha \eta_\alpha^{p - 1} \nabla_g u_\alpha \cdot \nabla_g \eta_\alpha\; dv_g \right|
\leq \varepsilon \int_M |\nabla_g u_\alpha|^p \eta_\alpha^p\; dv_g + c_\varepsilon \int_M |\nabla_g \eta_\alpha|^p u_\alpha^p\; dv_g\, .
\]

\n Also, by (\ref{3dqp}) and (\ref{lim}), it follows that

\begin{equation}\label{3dm}
A_\alpha \int_M|\nabla_g \eta_\alpha|^p u_\alpha^p\; dv_g \leq
A_\alpha (A_\alpha^{-\frac{1}{p}}
u_\alpha(y_\alpha)^{\frac{r-p}{p}})^p \int_{D_{\alpha}}
u_\alpha^p\; dv_g \leq c u_\alpha(y_\alpha)^r
(A_\alpha^{\frac{1}{p}} u_\alpha(y_\alpha)^{\frac{p-r}{p}})^n \leq
c m_\alpha^{- \varrho}\; .
\end{equation}

\n Consequently, combining these inequalities with (\ref{3dpt}),
one arrives at

\begin{equation}\label{31d}
A_\alpha \int_M |\nabla_g u_\alpha|^p \eta_\alpha^p\; dv_g + c
\alpha A_\alpha^{\frac{\tau}{p}} ||u_\alpha||_{L^p(M)}^{\tau - p}
\int_M u_\alpha^p \eta_\alpha^p\; dv_g + c \frac{\int_M u_\alpha^q
\eta_\alpha^p\; dv_g}{\int_M u_\alpha^q dv_g} \leq c
m_\alpha^{-\varrho}\, .
\end{equation}

\n Now, since $p > 1$,the non-sharp Riemannian Gagliardo-Nirenberg inequality produces

\begin{equation}\label{3ddb}
\left( \int_{B(y_\alpha, \frac{1}{2} A_\alpha^{\frac{1}{p}}
u_\alpha(y_\alpha)^{\frac{p-r}{p}})} u_\alpha^r\; dv_g
\right)^{\frac{p}{r \theta}} \leq \left( \int_M (u_\alpha
\eta_\alpha^p)^r\; dv_g \right)^{\frac{p}{r \theta}} \leq c
\left(\int_M |\nabla_g u_\alpha|^p \eta_\alpha^p\; dv_g\right)
\left( \int_M (u_\alpha \eta_\alpha^p)^q\; dv_g \right)^{\frac{p(1
- \theta)}{\theta q}}
\end{equation}

\[
+ c \left(\int_M |\nabla_g \eta_\alpha|^p u_\alpha^p\; dv_g\right) \left(\int_M (u_\alpha \eta_\alpha^p)^q\; dv_g \right)^{\frac{p(1 -
\theta)}{\theta q}} + c \left(\int_M (u_\alpha \eta_\alpha^p)^p\; dv_g\right) \left( \int_M (u_\alpha \eta_\alpha^p)^q\; dv_g
\right)^{\frac{p(1 - \theta)}{\theta q}}\, .
\]

\n Thanks to (\ref{3dm}), (\ref{31d}) and since $q \geq 1$, we then can estimate
each term of the right-hand side of (\ref{3ddb}). Indeed, we have

\[
\int_M |\nabla_g u_\alpha|^p \eta_\alpha^p \; dv_g \left( \int_M (u_\alpha \eta_\alpha^p)^q\; dv_g \right)^{\frac{p(1 -
\theta)}{\theta q}}
\leq A_\alpha \int_M |\nabla_g u_\alpha|^p \eta_\alpha^p\;
dv_g \left( \frac{\int_M u_\alpha^q \eta_\alpha^p\;
dv_g}{\int_M u_\alpha^q dv_g} \right)^{\frac{p(1 - \theta)}{\theta
q}} \leq c m_\alpha^{-\varrho(1 + \frac{p( 1 - \theta)}{\theta
q})},
\]

\[
\int_M |\nabla_g \eta_\alpha|^p u_\alpha^p\; dv_g \left(\int_M (u_\alpha \eta_\alpha^p)^q\; dv_g \right)^{\frac{p(1 -
\theta)}{\theta q}}
\leq A_\alpha \int_M |\nabla_g \eta_\alpha|^p u_\alpha^p\; dv_g
\left( \frac{\int_M u_\alpha^q \eta_\alpha^p\; dv_g}{\int_M
u_\alpha^q dv_g} \right)^{\frac{p(1 - \theta)}{\theta q}} \leq c
m_\alpha^{-\varrho(1 + \frac{p(1 - \theta)}{\theta q})} \; .
\]

\n By (\ref{dbd}) and because $p \geq \tau$, we have

\[
\alpha A_\alpha^{\frac{\tau - p}{p}} ||u_\alpha||_{L^p(M)}^{\tau -
p} > c \alpha^{\frac{p}{\tau}} > c > 0 \; ,
\]

\n so

\[
\int_M (u_\alpha \eta_\alpha^p)^p\; dv_g \left( \int_M (u_\alpha \eta_\alpha^p)^q\; dv_g \right)^{\frac{p(1 -
\theta)}{\theta q}}
\leq A_\alpha \int_M u_\alpha^p \eta_\alpha^p\; dv_g \left(\frac{
\int_M u_\alpha^q \eta_\alpha^p \; dv_g}{\int_M u_\alpha^q dv_g}
\right)^{\frac{p(1 - \theta)}{\theta q}} \leq c
m_\alpha^{-\varrho(1 + \frac{p(1 - \theta)}{\theta q})}\ .
\]

\n Replacing these three estimates in (\ref{3ddb}), one gets

\[
\left(\int_{B(y_\alpha, \frac{1}{2} A_\alpha^{\frac{1}{p}}
u_\alpha(y_\alpha)^{\frac{p-r}{p}})} u_\alpha^r\; dv_g
\right)^{\frac{p}{r \theta}} \leq c m_\alpha^{-\varrho(1 +
\frac{p(1 - \theta)}{\theta q})}\, ,
\]

\n so that

\[
m_\alpha^\varrho \int_{B(y_\alpha, \frac{1}{2}
A_\alpha^{\frac{1}{p}} u_\alpha(y_\alpha)^{\frac{p-r}{p}})}
u_\alpha^r\; dv_g \leq c m_\alpha^{\varrho(1 - \frac{r \theta}{p}
- \frac{r (1 - \theta)}{q})}\, .
\]

\n Since $m_\alpha \rightarrow \infty$ and

\[
1 - \frac{r \theta}{p} - \frac{r (1 - \theta)}{q} < 0\, ,
\]

\n we derive

\[
m_\alpha^\varrho \int_{B(y_\alpha, \frac{1}{2}
A_\alpha^{\frac{1}{p}} u_\alpha(y_\alpha)^{\frac{p-r}{p}})}
u_\alpha^r\; dv_g \rightarrow 0\, ,
\]

\n when $\alpha \rightarrow \infty$.  But this contradicts
(\ref{3dcm}).

\subsection{The final argument in the proof of Theorem \ref{tgno1}} \label{final-contradiction}

In the sequel, we will perform several estimates by using the
$L^r$-concentration and the pointwise estimation. By the scale invariance of the problem, we can assume that
 the radius of injectivity of $M$ grater than one.

Let $\eta \in C^1_0(\R)$ be a cutoff function as in the previous
section and define $\eta_{\alpha}(x) = \eta(d_g(x,x_\alpha))$.
From the inequality $GN_E(A(p,q,r,n))$ and by (\ref{desinf}), we have

\[
\left( \int_{B(0,1)} u_\alpha^r \eta_{\alpha}^r\; dx
\right)^{\frac{p}{r \theta}} \leq A(p,q,r,n) \left( \int_{B(0,1)}
|\nabla (u_\alpha \eta_{\alpha})|^p\; dx\right)
\left(\int_{B(0,1)} u_\alpha^q \eta_{\alpha}^q\; dx
\right)^{\frac{p(1 - \theta)}{\theta q}}
\]

\[
\leq A_{opt}^{\frac{p}{\tau}} \left( \int_{B(0,1)} |\nabla
(u_\alpha \eta_{\alpha})|^p\; dx\right) \left(\int_{B(0,1)}
u_\alpha^q \eta_{\alpha}^q\; dx \right)^{\frac{p(1 -
\theta)}{\theta q}}\, .
\]

\n Expanding the metric $g$ in normal coordinates around $x_\alpha$, one locally gets

\begin{equation} \label{car}
(1 - c d_g(x,x_\alpha)^2)\; dv_g \leq dx \leq (1 + c d_g(x,x_\alpha)^2)\; dv_g
\end{equation}

\n and

\begin{equation} \label{car1}
|\nabla(u_\alpha \eta_\alpha)|^p \leq |\nabla_g(u_\alpha \eta_\alpha)|^p(1 + c d_g(x,x_\alpha)^2) \, .
\end{equation}

\n Thus,

\[
\left( \int_{B(0,1)} u_\alpha^r \eta_{\alpha}^r\; dx
\right)^{\frac{p}{r \theta}}
\]

\[
\leq  \left( A_\alpha A_{opt}^{\frac{p}{\tau}} \int_M |\nabla_g
(u_\alpha \eta_{\alpha})|^p\; dv_g + c A_\alpha \int_M |\nabla_g
(u_\alpha\eta_{\alpha})|^p d_g(x,x_\alpha)^2\; dv_g \right)
\left(\frac{\int_{B(0,1)} u_\alpha^q \eta_{\alpha}^q\; dx}{\int_M
u_\alpha^q\; dv_g}\right)^{\frac{p(1 - \theta)}{\theta q}}\, .
\]

\n Applying then the inequalities

\[
|\nabla_g (u_\alpha \eta_{\alpha})|^p \leq |\nabla_g u_\alpha|^p
\eta_{\alpha}^p + c |\eta_{\alpha} \nabla_g u_\alpha|^{p - 1}
|u_\alpha \nabla_g \eta_{\alpha}| + c |u_\alpha \nabla_g
\eta_{\alpha}|^p \; ,
\]

\n (\ref{3ep}) and (\ref{rm}), we have

\[
A_{opt}^{\frac{p}{\tau}} \; A_\alpha \int_M |\nabla_g
u_\alpha|^p\; dv_g  \leq 1 - \alpha \left( A_\alpha \int_M
u_\alpha^p\; dv_g \right)^{\frac{\tau}{p}}\, ,
\]

\n one easily checks that

\[
\left( \int_{B(0,1)} u_\alpha^r \eta_{\alpha}^r\; dx
\right)^{\frac{p}{r \theta}}
\]
\begin{equation}\label{3du}
\leq  \left( 1 -  \alpha \left(A_\alpha \int_M u_\alpha^p\;
dv_g\right)^{\frac{\tau}{p}} + c F_\alpha + c G_\alpha + c
A_\alpha \int_{B(x_\alpha,1) \setminus B(x_\alpha,\frac{1}{2})}  u_\alpha^p \; dv_g \right)
\left(\frac{\int_{B(0,1)} u_\alpha^q \eta_{\alpha}^q\; dx}{\int_M
u_\alpha^q\; dv_g}\right)^{\frac{p(1 - \theta)}{\theta q}}\, ,
\end{equation}

\n where

\[
F_\alpha = A_\alpha \int_M |\nabla_g u_\alpha|^p \eta_{\alpha}^p
d_g(x,x_\alpha)^2\; dv_g
\]

\n and

\[
G_\alpha = A_\alpha \int_M |\nabla_g u_\alpha|^{p-1}
\eta_{\alpha}^{p-1} u_\alpha |\nabla_g \eta_{\alpha}|\; dv_g\, .
\]

\n We now estimate $F_\alpha$ and $G_\alpha$. Note that by
(\ref{eq1}), taking $u_\alpha$ as test function, we have

\[
A_\alpha \int_M |\nabla_gu_\alpha|^p dv_g \leq \lambda_\alpha \leq
A_{opt}^{-\frac{p}{\tau}} \; .
\]

\n Then applying H\"older inequality, (\ref{ic}), the definition of $\varphi_\alpha$,
the inequality above and the pointwise estimate, we obtain

\[
G_\alpha \leq \left(A_\alpha \int_M |\nabla_g u_\alpha|^p dv_g
\right)^{\frac{p - 1}{p}} \left( A_\alpha \int_{B(x_\alpha,1)
\setminus B(x_\alpha,\frac{1}{2})} u_\alpha^p dv_g
\right)^{\frac{1}{p}} \leq c \left( A_\alpha
\int_{B(x_\alpha,1) \setminus B(x_\alpha,\frac{1}{2})} u_\alpha^p
d_g(x,x_\alpha)^p dv_g \right)^{\frac{1}{p}}
\]

\begin{equation}\label{3dfq1}
\leq c \left(a_\alpha^{\frac{np - nr + pr}{r} + p
-\frac{np}{r} + n} \int_{B(0,a_\alpha^{-1}) \setminus B(0,
\frac{a_\alpha^{-1}}{2})} \varphi_\alpha^p |x|^p dh_\alpha
\right)^{\frac{1}{p}} \leq c_\lambda a_\alpha^2 \left(\int_{\R^n
\setminus B(0,1)} |x|^{p(1 - \lambda)} dx \right)^{\frac{1}{p}}
\leq c \; a_\alpha^2 \; ,
\end{equation}

\n with $\lambda$ being large enough. Similarely,

\[
A_\alpha \int_M |\nabla_g u_\alpha|^{p - 1} \eta_{\alpha}^p
u_\alpha d_g(x,x_\alpha)\; dv_g \leq \left(A_\alpha \int_M
|\nabla_g u_\alpha|^p dv_g \right)^{\frac{p - 1}{p}}
\left(A_\alpha \int_{B(x_\alpha,1)} u_\alpha^p d_g(x,x_\alpha)^p
dv_g\right)^{\frac{1}{p}}
\]
\begin{equation}\label{3dfq2}
\leq c a_\alpha^2 \left(\int_{B(0,a_\alpha^{-1})} \varphi_\alpha^p
|x|^p dh_\alpha \right)^{\frac{1}{p}} \leq c \; a_\alpha^2 \left(1
+ \int_{\R^n \setminus B(0,1)} |x|^{p(1 - \lambda)} dx
\right)^{\frac{1}{p}} \leq c \; a_\alpha^2\; .
\end{equation}

\n Now taking $u_\alpha d_g^2 \eta_\alpha^p$ as a test
function in (\ref{3ep}), one easily checks that

\[
F_\alpha = A_\alpha \int_M |\nabla_g u_\alpha|^p \eta_{\alpha}^p
d_g(x,x_\alpha)^2\; dv_g \leq c \int_{B(x_\alpha, 1)} u_\alpha^r
d_g(x,x_\alpha)^2\; dv_g + c A_\alpha \int_M |\nabla_g
u_\alpha|^{p - 1} \eta_{\alpha}^p u_\alpha d_g(x,x_\alpha)\; dv_g
+ c G_\alpha\, .
\]

\n Therefore, by (\ref{3dfq1}) and (\ref{3dfq2}),

\[
F_\alpha \leq c \int_{B(x_\alpha, 1)} u_\alpha^r
d_g(x,x_\alpha)^2\; dv_g + c a_\alpha^2 \, .
\]

\n Again using the pointwise estimate, we get

\begin{equation}\label{r1}
\int_{B(x_\alpha,1)} u_\alpha^r d_g(x,x_\alpha)^2\; dv_g \leq c
a_\alpha^ 2 \left(1 + \int_{B(0,a_\alpha^{-1}) \setminus B(0,1)}
\varphi_\alpha^r |x|^2 dx\right) \leq c a_\alpha^2 \left(1 + \int_{\R^n \setminus
B(0,1)} |x|^{2 - \lambda r} dx\right) \leq c a_\alpha^2
\end{equation}

\n for $\lambda$ big enough. Consequently,

\begin{equation}\label{z}
F_\alpha \leq c a_\alpha^2 \hspace{0,3cm} \mbox{and}
\hspace{0,3cm} G_\alpha \leq c a_\alpha^2 \, .
\end{equation}

\n Proceeding as in (\ref{3dfq1}) and since $p > 1$, we get

\[
A_\alpha \int_{B(x_\alpha,1)
\setminus B(x_\alpha,\frac{1}{2})} u_\alpha^p dv_g \leq c A_\alpha \int_{B(x_\alpha,1)
\setminus B(x_\alpha,\frac{1}{2})} u_\alpha^p d_g(x,x_\alpha)^p dv_g \leq c  a_\alpha^{2p} \leq c a_\alpha^2\; .
\]

\n Inserting this estimate and (\ref{z}) in (\ref{3du}), one arrives at

\begin{equation}\label{r3}
\left( \int_{B(x_\alpha,1)} u_\alpha^r \eta_\alpha^r\; dx
\right)^{\frac{p}{r \theta}} \leq  \left( 1 - \alpha
\left(A_\alpha \int_M u_\alpha^p\; dv_g\right)^{\frac{\tau}{p}} +
c a_\alpha^2 \right) \left(\frac{\int_{B(x_\alpha,1)} u_\alpha^q
\eta_\alpha^q\; dx}{\int_M u_\alpha^q\; dv_g}\right)^{\frac{p(1 -
\theta)}{\theta q}}\, .
\end{equation}

\n Now by (\ref{car}) and the mean value theorem, we
obtain

\[
\left(\int_M u_\alpha^r \eta_\alpha^r\; dx \right)^{\frac{p}{r \theta}} \geq \left( \int_M u_\alpha^r \eta_\alpha^r\; dv_g - c
\int_M u_\alpha^r \eta_\alpha^r d_g(x,x_\alpha)^2\; dv_g \right)^{\frac{p}{r \theta}}
\]

\[
\geq 1 - c\int_{M \setminus B(x_\alpha,1)} u_\alpha^r\; dv_g - c \int_M u_\alpha^r \eta_\alpha^r d_g(x,x_\alpha)^2\; dv_g
\]

\n and

\[
\left( \frac{\int_{B(x_\alpha,1)} u_\alpha^q \eta_\alpha^q\; dx}{\int_M u_\alpha^q\; dv_g} \right)^{\frac{p(1 - \theta)}{\theta q}} \leq \left(
\frac{\int_M u_\alpha^q \eta_\alpha^q\; dv_g + c \int_M u_\alpha^q \eta_\alpha^q d_g(x,x_\alpha)^2\; dv_g}{\int_M u_\alpha^q\; dv_g}
\right)^{\frac{p(1 - \theta)}{\theta q}}
\]

\[
\leq \left( \frac{\int_M u_\alpha^q \eta_\alpha^q\; dv_g}{\int_M u_\alpha^q\; dv_g} \right)^{\frac{p(1 - \theta)}{\theta q}} + c
\frac{\int_M u_\alpha^q \eta_\alpha^q d_g(x,x_\alpha)^2\; dv_g}{\int_M u_\alpha^q\; dv_g} \leq 1 + c \frac{\int_M u_\alpha^q
\eta_\alpha^q d_g(x,x_\alpha)^2\; dv_g}{\int_M u_\alpha^q\; dv_g}\, .
\]

\n Replacing these two estimates in (\ref{r3}), one gets

\[
\alpha \left(A_\alpha \int_M u_\alpha^p\;
dv_g\right)^{\frac{\tau}{p}} \leq c a_\alpha^2 + c \frac{\int_M
u_\alpha^q \eta_\alpha^q d_g(x,x_\alpha)^2\; dv_g}{\int_M
u_\alpha^q\; dv_g} + c \int_M u_\alpha^r \eta_\alpha^r
d_g(x,x_\alpha)^2\; dv_g + c \int_{M \setminus B(x_\alpha,1)}
u_\alpha^r\; dv_g\, .
\]

\n Using the pointwise estimate with $\lambda r - n = 2$, we have

\[
\int_{M \setminus B(x_\alpha,1)} u_\alpha^r\; dv_g \leq c \int_{M
\setminus B(x_\alpha,1)} u_\alpha^r d(x,x_\alpha)^{\lambda r} \;
dv_g \leq c_\lambda a_\alpha^{\lambda r - n} = c \, a_\alpha^2 \;.
\]

\n So, by this estimate and (\ref{r1}), one concludes that

\begin{equation}\label{eff}
\alpha \left(A_\alpha \int_M u_\alpha^p\;
dv_g\right)^{\frac{\tau}{p}} \leq c a_\alpha^2 + c \frac{\int_M
u_\alpha^q \eta_\alpha^q d_g(x,x_\alpha)^2\; dv_g}{\int_M
u_\alpha^q\; dv_g} \; .
\end{equation}

\n By the pointwise estimate, for $\lambda$ big enough

\[
\frac{\int_M u_\alpha^q \eta_\alpha^q d_g(x,x_\alpha)^2\;
dv_g}{\int_M u_\alpha^q\; dv_g} \leq \frac{\int_{B(x_\alpha,1)}
u_\alpha^q d_g(x,x_\alpha)^2\; dv_g}{\int_M u_\alpha^q\; dv_g}
\]
\[= a_\alpha^2 \int_{B(0,a_\alpha^{-1})} \varphi^q_\alpha |x|^2
dh_\alpha \leq c a_\alpha^2 \left( 1 + c_\lambda \int_{\R^n
\setminus B(0,1)} |x|^{2 - \lambda q} dx\right) \leq c \;
a_\alpha^2 \; .
\]

\n Introducing this inequality in (\ref{eff}), we readily deduce
that

\[
\alpha \left(A_\alpha \int_M u_\alpha^p\;
dv_g\right)^{\frac{\tau}{p}} \leq c \; a_\alpha^2 \, .
\]

\n Then, since

\[
\int_M u_\alpha^p \; dv_g \geq \int_{B(x_\alpha,a_\alpha)} u_\alpha^p\; dv_g = a_\alpha^{\frac {n(r - p)}{r}} \int_{B(0,1)}
\varphi_\alpha^p\; dh_\alpha \geq c a_\alpha^{\frac {n(r - p)}{r}}
\]

\n for $\alpha$ large enough and since $a_\alpha^{\frac{n(r - p)}{r}} A_\alpha = a_\alpha^p$, one arrives at

\[
\alpha a_\alpha^\tau \leq c \; a_\alpha^2 \, .
\]

\n Finally, because $\tau \leq 2$, we arrive at
the contradiction

\[
\alpha \leq c\; a_\alpha^{2 - \tau} \; .
\]

\bl

\section{Existence of extremals for the optimal Gagliardo-Nirenberg inequality}

Let us now prove Theorem \ref{extremal}

By Theorem \ref{tgno1} we have that the inequality

\[
\left( \int_M |u|^r\; dv_g \right)^{\frac{\tau}{r \theta}} \leq
\left( A_{opt} \left(\int_M |\nabla_g u|^p\;
dv_g\right)^{\frac{\tau}{p}} + B_{opt} \left(\int_M |u|^p\;
dv_g\right)^{\frac{\tau}{p}} \right) \left( \int_M |u|^q\; dv_g
\right)^{\frac{\tau(1 - \theta)}{\theta q}}\, ,
\]

\n is valid for all $u \in H^{1,p}(M)$.

Let $\alpha > 0$ and $c_\alpha = B_{opt} - \alpha^{-1}$. Define

\[
J_\alpha(u) = \left( \int_M |u|^r\; dv_g\right)^{\frac{\tau}{r
\theta}}\left(\int_M |u|^q dv_g \right)^{- \frac{\tau(1 -
\theta)}{\theta q}} - c_\alpha \left(\int_M |u|^p\; dv_g
\right)^{\frac{\tau}{p}} \, .
\]

\n By the definition of $B_{opt}$ we have

\[
\nu_\alpha = \sup_{u \in E} J_\alpha (u) > A_{opt}\, ,
\]

\n where $E = \{ u \in H^{1,p}(M):\; ||\nabla_g u||_{L^p(M)} = 1
\}$. Clearly, this supremum is well-defined.

In the same fashion as in section \ref{EL}, we will give the proof for the
$C^1$-case $q>1$, the case $q=1$ being taken care or by the methods contained in
\cite{Hu}. Since then $J_\alpha$ is of class $C^1$, by
using standard variational arguments, we find a maximizer
$\tilde{u}_\alpha \in E$ of $J_\alpha$, and then

\[
J_\alpha(\tilde{u}_\alpha) = \nu_\alpha = \sup_{u \in E}
J_\alpha(u)\ .
\]

\n The function $\tilde{u}_\alpha$ satisfies the Euler-Lagrange
equation

\[
\frac{1}{\theta}||\tilde{u}_\alpha||_{L^r(M)}^{\frac{\tau - r
\theta}{\theta}} ||\tilde{u}_\alpha||_{L^q(M)}^{- \frac{\tau(1 -
\theta)}{\theta}} \tilde{u}_\alpha^{r - 1} - \frac{(1 -
\theta)}{\theta}
||\tilde{u}_\alpha||_{L^r(M)}^{\frac{\tau}{\theta}}
||\tilde{u}_\alpha||_{L^q(M)}^{-\frac{\tau(1 + \theta) + \theta
q}{\theta}} \tilde{u}_\alpha^{q - 1} - c_\alpha
||\tilde{u}_\alpha||_{L^p(M)}^{\tau - p} \tilde{u}_\alpha^{p - 1}
= \nu_\alpha \Delta_{p,g} \tilde{u}_\alpha \; ,
\]

\n where $\Delta_{p,g} = -{\rm div}_g(|\nabla_g|^{p-2} \nabla_g)$
is the $p$-Laplace operator of $g$. Because $\nabla_g
|\tilde{u}_\alpha| = \pm \nabla_g \tilde{u}_\alpha$, we can assume
$\tilde{u}_\alpha \geq 0$ and again be  Tolksdorf regularity \cite{To}, we have  that
 $\tilde{u}_\alpha \in C^1(M)$. Setting $u_\alpha =
\frac{\tilde{u}_\alpha}{||\tilde{u}_\alpha||_{L^r(M)}}$,
we find

\begin{equation} \label{ext1}
\lambda_\alpha^{-1} A_\alpha \Delta_{p,g} u_\alpha + c_\alpha
 A_\alpha^{\frac{\tau}{p}} ||u_\alpha||_{L^p(M)}^{\tau -
p} u_\alpha^{p - 1} + \frac{1 - \theta}{\theta}
||u_\alpha||_{L^q(M)}^{-q} u_\alpha^{q - 1} = \frac{1}{\theta}
u_\alpha^{r - 1}\ \ {\rm on}\ \ M\, ,
\end{equation}

\n where  $||u_\alpha||_{L^r(M)} = 1$,

\[
A_\alpha = \left(\int_M u_\alpha^q\; dv_g \right)^{\frac{p(1 -
\theta)}{\theta q}}
\]

\n and

\[
\lambda_\alpha = \nu_\alpha^{-1} ||u_\alpha||_{L^q(M)}^{\frac{(p -
\tau)(1 - \theta)}{\theta}} ||\tilde{u}_\alpha||_{L^r(M)}^{\tau -
p} \; .
\]

\n Up to taking a subsequence, we can assume that there exists $A \in \R$ such that

\[
\lim_{\alpha \rightarrow \infty} A_\alpha = A \; .
\]

\n Then there are two possibilities: either $A=0$ or $A>0$. We presently show that $A=0$ cannot happen. Indeed,
if $A=0$ then (\ref{3}) of Theorem \ref{tgno1} holds. Then we can follow the proof of theorem \ref{tgno1} and we get
\[
c_\alpha \leq c \; a_\alpha^{2 - \tau} \; ,
\]

\n as in the end of section \ref{final-contradiction}. Since $\tau < 2$, $c_\alpha = B_{opt}
- \alpha^{-1}$ and $\lim_{\alpha \rightarrow \infty} a_\alpha = 0$, we have that above inequality contradicts (\ref{sco}). Therefore $A>0$.

\n Using this fact in
 (\ref{ext1}), we see that there exists $c > 0$ such that

\[
\int_M|\nabla_g u_\alpha|^p dv_g + \left(\int_M u_\alpha^p
dv_g\right)^{\frac{\tau}{p}} \leq c \;
\]

\n for all $\alpha$. Then, up to a subsequence, $u_\alpha
\rightharpoonup u_0$ in $H^{1,p}(M)$. Since $||u_\alpha||_{L^r(M)}
= 1$ for all $\alpha$, we have that $||u_0||_{L^r(M)} = 1$.

\n From  (\ref{ext1}) we have

\[
\int_M |\nabla_g u_\alpha|^{p - 2} \nabla_g u_\alpha \nabla_g h
dv_g \leq c \int_M u_\alpha^{r - 1} h dv_g \; ,
\]

\n for an arbitrary test function $h \geq 0$. Then by Moser's iterative scheme,

\[
\sup_{x \in M} u_\alpha \leq c \left(\int_M u_\alpha^r dv_g
\right)^{\frac{1}{r}} \leq c \; ,
\]

\n for all $\alpha$.

\n Using Tolksdorf regularity in (\ref{ext1}), it follows that $u_\alpha \rightarrow
u_0$ in $C^1(M)$.

The function $\tilde{u}_\alpha$ satisfies

\[
\left(\int_M \tilde{u}_\alpha^r dv_g \right)^{\frac{\tau}{r
\theta}} \geq \left( A_{opt}^{\frac{p}{\tau}} \left(\int_M
|\nabla_g \tilde{u}_\alpha|^p dv_g \right)^{\frac{\tau}{p}} +
(B_{opt} - \frac{1}{\alpha}) \left(\int_M \tilde{u}_\alpha^p dv_g
\right)^{\frac{\tau}{p}} \right) \left(\int_M \tilde{u}_\alpha^q
dv_g \right)^{\frac{\tau(1 - \theta)}{\theta q}} \; ,
\]

\n and since $u_\alpha =
\frac{\tilde{u}_\alpha}{||\tilde{u}_\alpha||_{L^r(M)}}$, we get

\[
1 \geq \left( A_{opt}^{\frac{p}{\tau}} \left(\int_M |\nabla_g
u_\alpha|^p dv_g \right)^{\frac{\tau}{p}} + (B_{opt} -
\frac{1}{\alpha}) \left(\int_M u_\alpha^p dv_g
\right)^{\frac{\tau}{p}} \right) \left(\int_M u_\alpha^q dv_g
\right)^{\frac{\tau(1 - \theta)}{\theta q}} \; .
\]

\n Taking the limit in $\alpha$ in this inequality, we find

\[
1 \geq \left( A_{opt}^{\frac{p}{\tau}} \left(\int_M |\nabla_g
u_0|^p dv_g \right)^{\frac{\tau}{p}} + B_{opt} \left(\int_M u_0^p
dv_g \right)^{\frac{\tau}{p}} \right) \left(\int_M u_0^q dv_g
\right)^{\frac{\tau(1 - \theta)}{\theta q}} \; .
\]

\n And then $u_0$ is an extremal function for
$GN_R(A_{opt},B_{opt})$. \bl

\end{document}